\documentclass[11pt]{amsart}
\usepackage[left=2.3cm, top=2.3cm,bottom=2.3cm,right=2.3cm]{geometry}

\usepackage[T1]{fontenc}
\usepackage{lmodern}
\usepackage{amssymb}
\usepackage{stmaryrd}
\usepackage{microtype}
\usepackage{booktabs}
\usepackage{xcolor}
\usepackage[hidelinks]{hyperref}

\newtheorem{theorem}{Theorem}[section]
\newtheorem{proposition}[theorem]{Proposition}
\newtheorem{lemma}[theorem]{Lemma}
\newtheorem{corollary}[theorem]{Corollary}

\theoremstyle{definition}

\theoremstyle{remark}
\newtheorem{remark}[theorem]{Remark}

\newcommand{\F}{\mathcal L}
\newcommand{\T}{\mathcal T}
\newcommand{\RR}{\mathbb R}
\newcommand{\CC}{\mathbb C}
\newcommand{\LP}{\mathcal{LP}_{+}}
\newcommand{\lsem}{\llbracket}
\newcommand{\rsem}{\rrbracket}

\newcommand{\parhead}[1]{%
  \par\addvspace{0.5\baselineskip}%
  \noindent\textbf{#1}\enspace\ignorespaces}

\newcommand{\doi}[1]{\href{https://doi.org/#1}{doi:\,#1}}

\title[Real zeros of binomial-power polynomials and the Le Roy function]{%
Real-rootedness of $\sum_{k=0}^n\binom nk^\gamma z^k$
and $\sum_{k=0}^{\infty}\frac{z^k}{(k!)^\gamma}$}

\author{Zakhar Kabluchko}
\address{Zakhar Kabluchko: Institut f\"ur Mathematische Stochastik,
Universit\"at M\"unster,
Orl\'eans-Ring 10,
48149 M\"unster, Germany}
\email{zakhar.kabluchko@uni-muenster.de}

\author{Alexander Marynych}
\address{Alexander Marynych: Department of Mathematics, Kyiv School of Economics, Kyiv, Ukraine}
\email{omarynych@kse.org.ua}

\date{\today}



\hypersetup{
  pdftitle={Real-rootedness of coefficientwise powers of binomial polynomials and the Le Roy function}
}

\begin{document}

\begin{abstract}
We prove that for every real \(\gamma\ge2\) and every integer
\(n\ge1\), the polynomial
\[
 P_{n,\gamma}(z)=\sum_{k=0}^n\binom nk^\gamma z^k
\]
has \(n\) simple negative zeros. Moreover, the zeros of
\(P_{n,\gamma}\) and \(P_{n+1,\gamma}\) strictly interlace.
In contrast, for each fixed \(\gamma \in (0,2)\backslash\{1\}\), the number of
nonreal zeros of \(P_{n,\gamma}\) tends to infinity as
\(n\to\infty\).
We also show that for \(\gamma\ge2\), the  Le Roy
function
\[
 \F_\gamma(z)=\sum_{k=0}^{\infty}\frac{z^k}{(k!)^\gamma}
\]
has infinitely many zeros, all simple and negative.  On the other hand, for \(\gamma \in (0,2)\backslash\{1\}\) the function $\F_\gamma$ has infinitely many nonreal zeros.   Furthermore, we show that for a real number \(s\), the sequence
\((1/(k!)^s)_{k\ge0}\) is a multiplier sequence if
and only if \(s\in\{0\}\cup[1,\infty)\).
\end{abstract}

\maketitle

\section{Main results and related work}\label{sec:intro}

\subsection{Main results}
For every integer \(n\ge1\) and every real \(\gamma>0\), let \(P_{n,\gamma}\) be a polynomial obtained by raising each coefficient of \((1+x)^n\) to the power \(\gamma\), that is
\begin{equation}\label{eq:objects}
 P_{n,\gamma}(x)=\sum_{k=0}^n\binom nk^\gamma x^k.
\end{equation}
A polynomial is called
\emph{real-rooted} if all its zeros are real. In this paper, we study real-rootedness of the polynomials $P_{n,\gamma}$ and related families of polynomials and entire functions.

\begin{theorem}[Zeros of $P_{n,\gamma}$]\label{th:main}\label{th:sharp}
Let \(\gamma>0\) be real.
\begin{enumerate}
\item If \(\gamma\ge2\), then for every integer \(n\ge1\),
the polynomial \(P_{n,\gamma}\) has \(n\) simple negative zeros.
Writing these zeros in increasing order as
\(
 x_{n,1}<\cdots<x_{n,n}<0,
\)
the zeros of \(P_{n,\gamma}\) and \(P_{n+1,\gamma}\) strictly
interlace:
\[
 x_{n+1,1}<x_{n,1}<x_{n+1,2}<\cdots
 <x_{n,n}<x_{n+1,n+1}<0.
\]
\item If \(1<\gamma<2\), then the number of nonreal zeros of
\(P_{n,\gamma}\), counted with multiplicity, tends to infinity
as \(n\to\infty\).
\item If \(\gamma=1\), then \(P_{n,1}(x)=(1+x)^n\) for every
integer \(n\ge1\).
\item If \(0<\gamma<1\), then \(P_{n,\gamma}\) has no real
zeros when \(n\) is even and exactly one real zero when \(n\)
is odd. In the latter case, the zero is \(-1\) and is simple.
\end{enumerate}
In particular, \(P_{n,\gamma}\) is real-rooted for every integer
\(n\ge1\) if and only if $\gamma\in\{1\}\cup[2,\infty)$.
\end{theorem}

For every real \(\gamma>0\), the \emph{Le Roy function} \cite{LeRoy,Gerhold} is the entire function defined by
\begin{equation}\label{eq:LeRoy-definition}
 \F_\gamma(z)=\sum_{k=0}^\infty\frac{z^k}{(k!)^\gamma}.
\end{equation}
Its connection
with the polynomials $P_{n,\gamma}$ is given by
\[
 P_{n,\gamma}\!\left(\frac{z}{n^\gamma}\right)
 \longrightarrow\F_\gamma(z)
 \qquad\text{as }n\to\infty,
\]
locally uniformly on \(\CC\), for every real \(\gamma>0\).

\begin{theorem}[Zeros of the Le Roy function]\label{th:LeRoy}
Let \(\gamma>0\) be real.
\begin{enumerate}
\item If \(\gamma\ge2\), then \(\F_\gamma\) has infinitely
many zeros, all simple and negative.
\item If \(1<\gamma<2\), then \(\F_\gamma\) has infinitely
many nonreal zeros and finitely many real zeros, all negative.
The total multiplicity of its real zeros is positive and odd.
\item If \(\gamma=1\), then \(\F_1(z)=e^z\), which has no zeros.
\item If \(0<\gamma<1\), then \(\F_\gamma\) has infinitely
many zeros, all nonreal.
\end{enumerate}
Moreover,
\[
\begin{aligned}
 \F_\gamma\in\LP
 &\quad\Longleftrightarrow\quad \F_\gamma\in\mathcal{LP}\\
 &\quad\Longleftrightarrow\quad
   \F_\gamma\text{ has no nonreal zeros}\\
 &\quad\Longleftrightarrow\quad
   \gamma\in\{1\}\cup[2,\infty).
\end{aligned}
\]
Here \(\mathcal{LP}\) denotes the Laguerre--P\'olya class,
consisting of locally uniform limits on \(\CC\) of real-rooted
real polynomials, and \(\LP\) denotes the class of locally uniform
limits on \(\CC\) of polynomials with nonnegative coefficients
and only nonpositive zeros.
\end{theorem}

A related family is given by the \emph{modified Laguerre polynomials}
\[
 Q_{n,s}(x)=\sum_{k=0}^n(-1)^k\binom nk\frac{x^k}{(k!)^s},
 \qquad n\ge1,\quad s\in\RR.
\]
At \(s=1\), we recover the ordinary
Laguerre polynomial \(Q_{n,1}=L_n\).

\begin{theorem}[Zeros of $Q_{n,s}$]\label{th:Q}
Let \(s\in\RR\).
\begin{enumerate}
\item If \(s\ge1\), then for every integer \(n\ge1\), the
polynomial \(Q_{n,s}\) has \(n\) simple positive zeros.
Writing these zeros in increasing order as
\(
 0<r_{n,1}<\cdots<r_{n,n},
\)
the zeros of \(Q_{n,s}\) and \(Q_{n+1,s}\) strictly interlace:
\[
 0<r_{n+1,1}<r_{n,1}<r_{n+1,2}<\cdots
 <r_{n,n}<r_{n+1,n+1}.
\]
\item If \(0<s<1\), then the number of nonreal zeros of
\(Q_{n,s}\), counted with multiplicity, tends to infinity
as \(n\to\infty\).
\item If \(s=0\), then \(Q_{n,0}(x)=(1-x)^n\) for every
integer \(n\ge1\).
\item If \(s<0\), then \(Q_{n,s}\) has no real zeros when
\(n\) is even and exactly one real zero when \(n\) is odd.
In the latter case, the zero is positive and simple.
\end{enumerate}
In particular, \(Q_{n,s}\) is real-rooted for every integer
\(n\ge1\) if and only if $s\in\{0\}\cup[1,\infty)$.
\end{theorem}

A real sequence
\((\lambda_k)_{k\ge0}\) is a \emph{multiplier sequence} if, for
every real-rooted real polynomial \(\sum_{k=0}^n a_kx^k\), the
polynomial \(\sum_{k=0}^n\lambda_ka_kx^k\) is real-rooted
or identically zero.

\begin{theorem}[Powers of the reciprocal factorial sequence]\label{th:multiplier}
For a real number \(s\), the sequence
$
 \left(1/(k!)^s\right)_{k\ge0}
$
is a multiplier sequence if and only if
\(
 s\in\{0\}\cup[1,\infty).
\)
\end{theorem}

A sequence
\((a_k)_{k\ge0}\) of nonnegative real numbers is called a
\emph{P\'olya frequency sequence} if the infinite Toeplitz
matrix \((a_{j-i})_{i,j\ge0}\), with \(a_k=0\) for \(k<0\),
is totally nonnegative, that is, all its minors are nonnegative.

\begin{corollary}
\label{cor:PF}
For a real number \(\gamma\), the following conditions
are equivalent:
\begin{enumerate}
\item The sequence
\(
 \left(1/(k!)^\gamma\right)_{k\ge0}
\)
is a P\'olya frequency sequence.
\item \(\gamma\in\{0,1\}\cup[2,\infty)\).
\end{enumerate}
\end{corollary}

We also determine the limiting empirical zero distributions of the polynomials
\(P_{n,\gamma}\) for \(\gamma\ge2\) and of
\(Q_{n,s}(n^s\,\cdot)\) for \(s\ge1\), as $n\to\infty$. In the latter case, the limits are given by the
Fuss--Catalan distribution. More generally, we show that the suitably normalized
multiplier operators associated with the sequences
\(((k!)^{-s})_{k\ge0}\), \(s\ge1\), transform limiting
zero distributions by free multiplicative convolution
with a Fuss--Catalan distribution. These results are given in  Section~\ref{sec:positions}.

\subsection{Historical background and related work}
\label{subsec:history}

Coefficientwise powers of real-rooted polynomials have been studied
in a more general setting. Fisk~\cite[Question~10, p.~722]{Fisk}
asked whether, for a real-rooted polynomial
\(f(x)=\sum_{k=0}^n a_kx^k\) with positive coefficients, the
polynomial \(\sum_{k=0}^n a_k^\gamma x^k\) is real-rooted for
every real \(\gamma\ge1\).
Wang and Zhang~\cite{WZ} proved this in degree three, gave a
counterexample in degree four, and showed that all sufficiently
large powers preserve real-rootedness in each fixed degree,
uniformly over the input polynomials.
Khare~\cite[Section~38.4, Question~38.4]{Khare} formulates the
problem of classifying the admissible powers and determining the
sharp degree-dependent thresholds.

Jalowy, Kabluchko, and
Marynych~\cite[Remark~4.8]{JKM} conjectured
real-rootedness of $P_{n,\gamma}$ for every integer \(n\ge1\) and every real
\(\gamma\ge2\), and reported numerical evidence of nonreal zeros
in sufficiently large degrees when \(1<\gamma<2\).
Theorem~\ref{th:main} confirms this conjecture.

The integer-parameter cases of the Le Roy function belong to the
classical theory of generalized hypergeometric functions.
A theorem of Hurwitz~\cite{Hurwitz}, independently rediscovered by
Hille~\cite{Hille}, states that
\({}_0F_q(\,;b_1,\ldots,b_q;z)\) belongs to \(\LP\) for every
integer \(q\ge0\) and all positive lower parameters
\(b_1,\ldots,b_q\); see Sokal~\cite[Theorem~1]{Sokal}.
Since
\[
 \F_m(z)={}_0F_{m-1}(\,;\underbrace{1,\ldots,1}_{m-1};z),
 \qquad m=1,2,\ldots,
\]
this gives \(\F_m\in\LP\) for every positive integer \(m\).
For this factorial family, Theorem~\ref{th:LeRoy} determines
the admissible positive real exponents and proves simplicity
of all zeros for \(\gamma\ge2\).

The asymptotic behavior of the Le Roy function was studied by
Le Roy~\cite{LeRoy} and Gerhold~\cite{Gerhold}.
Boudabsa and Simon~\cite[Section~5]{BS} studied its probabilistic
representations, the complete monotonicity of
\(x\mapsto\F_\gamma(-x)\) on \((0,\infty)\) for
\(0<\gamma\le1\), and its negative-axis asymptotics.
Coefficient criteria for real zeros of entire functions and their
Taylor sections have also been studied extensively.
For an entire function \(f(z)=\sum_{k=0}^\infty a_kz^k\)
with \(a_k>0\), Hutchinson~\cite{Hutchinson} proved that
\(f\) and all its Taylor sections have only real negative
zeros if the second quotients satisfy
\[
 q_k(f):=\frac{a_{k-1}^2}{a_{k-2}a_k}\ge4,
 \qquad k\ge2.
\]
Related coefficient criteria were developed by Katkova, Lobova,
and Vishnyakova~\cite{KLV} and by Nguyen and
Vishnyakova~\cite{NV2018,NV2024}.
For the Le Roy function,
\[
 q_k(\F_\gamma)=\left(\frac{k}{k-1}\right)^\gamma
 \longrightarrow1
 \qquad\text{as }k\to\infty.
\]
Consequently, the sufficient conditions requiring a fixed lower
bound greater than one on all second quotients do not apply
to this family.

The P\'olya--Schur theorem~\cite{PS} characterizes multiplier
sequences through their exponential generating functions,
whereas the Aissen--Schoenberg--Whitney--Edrei
theorem~\cite{ASW,Edrei} characterizes P\'olya frequency sequences
through their ordinary generating functions.
Multiplier sequences defined through ratios of Gamma functions
and other special functions were studied by Craven and
Csordas~\cite{CravenCsordas}. Katkova and
Vishnyakova~\cite{KV2024} investigated the related problem of
characterizing diagonal transformations that preserve the class
of P\'olya frequency sequences.
In a MathOverflow post of August 2024, motivated by a positivity
question for symmetric functions, Stanley~\cite{StanleyMO}
asked whether \(\bigl((k!)^{-m}\bigr)_{k\ge0}\) is a P\'olya
frequency sequence for every positive integer \(m\).
Kwa\'snicki's answer in the same discussion recalled the
Hurwitz--Hille theorem and an alternative argument using
P\'olya--Schur multiplier sequences.
Corollary~\ref{cor:PF} gives the corresponding classification
for all \emph{real} exponents.

The polynomials $P_{n,\gamma}$ are also related to the Brenke polynomials
studied by Dur\'an~\cite{Duran}. Indeed,
\[
 \F_\gamma(t)\F_\gamma(xt)
 =1+\sum_{n=1}^\infty
   \frac{P_{n,\gamma}(x)}{(n!)^\gamma}t^n,
\]
so their normalized versions form a Brenke family.
Dur\'an's results connect real-rootedness of these polynomial
families with Laguerre--P\'olya membership of their generating
functions.

\section{An inequality that detects a collision of zeros}\label{sec:dilation}
Consider the Euler differential operator \(D=x\,d/dx\). Because \(Dx^k=kx^k\), it is natural to define real powers of $D$ by
\begin{equation}\label{eq:D}
 D^s\Big(\sum_{k\ge0}c_kx^k\Big)=\sum_{k\ge1}k^s c_kx^k,
 \qquad s>0.
\end{equation}
This definition works for polynomials and entire functions. Multiplication
of their coefficients by \(k^s\) preserves an infinite radius of convergence.

\begin{lemma}[Key inequality]\label{lem:dil}
Let \(v(x)\) be a polynomial with real
coefficients or a real entire function (an entire function real-valued
on the real axis). If \(v(0)=v(r)=0\) for some \(r>0\), then for every \(1\le s\le2\) we have
\begin{equation}\label{eq:energy}
 \int_0^r v(x)D^s v(x)\,\frac{{\rm d}x}{x}\le0.
\end{equation}
\end{lemma}

\begin{proof}
First take \(v(x)\) to be a polynomial and \(1<s<2\). The number
\[
 C_s=\int_0^\infty (e^{-h}-1+h)\,\frac{{\rm d}h}{h^{1+s}}
\]
is positive and finite. Indeed, the expression $e^{-h}-1+h$ is positive
for \(h>0\), is \(O(h^2)\) at zero, and is \(O(h)\) at infinity. (Integration by parts twice yields \(C_s=\Gamma(-s)\), but we don't need this value for the proof.)
For every integer \(k\ge1\), the change of variables \(u=kh\) gives
\begin{equation}\label{eq:scaling-identity}
 \int_0^\infty (e^{-kh}-1+kh)\,\frac{{\rm d}h}{h^{1+s}}
 =k^s\int_0^\infty (e^{-u}-1+u)\,\frac{{\rm d}u}{u^{1+s}}
 =k^s C_s.
\end{equation}
Applying this identity to each monomial of \(v\) and
using~\eqref{eq:D}, we obtain
\begin{equation}\label{eq:dilation}
 C_sD^s v(x)=\int_0^\infty
 \{v(e^{-h}x)-v(x)+hDv(x)\}\,\frac{{\rm d}h}{h^{1+s}}.
\end{equation}

Define the inner product \(\langle f,g\rangle_r=\int_0^r fg\,{\rm d}x/x\), which is finite whenever polynomials (or entire functions) $f$ and $g$ vanish at $0$. Observe that, for $h>0$,
\[
 \|v(e^{-h}\,\cdot)\|_r^2
   =\int_0^{e^{-h}r}v(u)^2\,\frac{{\rm d}u}{u}\le\|v\|_r^2,
 \qquad
 \langle v,Dv\rangle_r=\int_0^r vv'\,{\rm d}x=0,
\]
where in the first formula we have used the scale invariance of the weight ${\rm d}x/x$, and in the second both boundary conditions $v(0)=v(r)=0$. Hence, by the Cauchy--Schwarz inequality,
\begin{equation}\label{eq:cauchy-schwarz-ineq}
 \langle v(e^{-h}\,\cdot)-v,v\rangle_r
 =\langle v(e^{-h}\,\cdot),v\rangle_r-\|v\|_r^2
 \leq\|v(e^{-h}\,\cdot)\|_r\|v\|_r-\|v\|_r^2
 \leq0.
\end{equation}
Taking the inner product of both sides of~\eqref{eq:dilation}
with \(v\), we obtain
\[
\begin{aligned}
 C_s\langle D^sv,v\rangle_r
 &=\int_0^r v(x)\int_0^\infty
 \bigl(v(e^{-h}x)-v(x)+hDv(x)\bigr)
 \frac{{\rm d}h}{h^{1+s}}\frac{{\rm d}x}{x}\\
 &=\int_0^\infty
 \bigl(\langle v(e^{-h}\,\cdot)-v,v\rangle_r
 +h\langle Dv,v\rangle_r\bigr)
 \frac{{\rm d}h}{h^{1+s}}\\
 &=\int_0^\infty
 \langle v(e^{-h}\,\cdot)-v,v\rangle_r
 \frac{{\rm d}h}{h^{1+s}}\\
 &\leq0.
\end{aligned}
\]
The interchange of integrals in the second equality is justified
by Fubini's theorem: expanding \(v\) into monomials, we see that
each term is absolutely integrable by the finiteness of \(C_s\)
and~\eqref{eq:scaling-identity}.
In the third equality we have used \(\langle Dv,v\rangle_r=0\),
and the last inequality follows from~\eqref{eq:cauchy-schwarz-ineq}.
Since \(C_s>0\), this proves~\eqref{eq:energy} for \(1<s<2\) when $v(x)$ is a polynomial.

\medskip
The inequality extends to the endpoints
\(s=1\) and \(s=2\) by continuity. Indeed, since \(v\) is a polynomial, it follows from~\eqref{eq:D}
that the map
$
 s\longmapsto
 \int_0^r v(x)D^sv(x)\,\frac{{\rm d}x}{x}
$
is continuous.  This proves the lemma for
polynomials.

\medskip
Let now \(v\) be a real entire function satisfying
\(v(0)=v(r)=0\). Write $v(x)=\sum_{k=1}^\infty c_kx^k$
and let
$
v_N(x)=\sum_{k=1}^N c_kx^k
$
be the \(N\)-th partial sum of its Taylor series at zero.
Set
\[
 w_N(x)=v_N(x)-\frac{v_N(r)}{r}x.
\]
By construction, \(w_N(0)=w_N(r)=0\).
Since the power series defining \(v(x)/x\) and \(D^sv(x)\)
converge locally uniformly and \(v_N(r)\to v(r)=0\),
we have
\[
 \frac{w_N(x)}{x}\to\frac{v(x)}{x},
 \qquad D^sw_N(x)\to D^sv(x),
\]
both uniformly on \([0,r]\) as \(N\to\infty\). The
quotients at zero are understood by continuity.
Applying~\eqref{eq:energy}, already proved for polynomials,
to \(w_N\), we obtain
\[
 \int_0^r w_N(x)D^sw_N(x)\,\frac{{\rm d}x}{x}\le0.
\]
Passing to the limit as \(N\to\infty\) using the above uniform
convergence yields~\eqref{eq:energy} for \(v\).
\end{proof}

\section{Jensen polynomials of the Le Roy function}
\label{sec:lag}
For every integer \(n\ge1\) and every real $s$,
consider  the polynomial
\begin{equation}\label{eq:Q}
 Q_{n,s}(x)=\sum_{k=0}^n(-1)^k\binom nk\frac{x^k}{(k!)^s}.
\end{equation}
For \(s>-1\), the polynomials \(Q_{n,s}(-x)\) are the
Jensen polynomials of the Le Roy function \(\F_{s+1}\);
see Section~\ref{sec:transfer}.

\medskip
Observe that for  \(s=1\) we have $Q_{n,1}(x)=\sum_{k=0}^n(-1)^k\binom nk\frac{x^k}{k!} = L_n(x)$
where \(L_n\) is the ordinary Laguerre polynomial of degree \(n\).
It is well known that \(L_n\) has \(n\) simple positive zeros;
see~\cite[Chapter V]{SzegoOP}.

\begin{proposition}\label{prop:Q}
For every integer \(n\ge1\) and every real \(s\in[1,2]\),
the polynomial \(Q_{n,s}\) has \(n\) simple positive zeros.
\end{proposition}

\begin{proof}
Fix \(n\) and abbreviate \(Q=Q_{n,s}\). Let \(c_k = (-1)^k\binom nk \frac{1}{(k!)^s}\) be
the coefficient of \(x^k\) in \(Q\). Then
\[
 k^{s+1}c_k=(k-1-n)c_{k-1},\qquad 1\le k\le n .
\]
Consider the polynomial \(v=DQ=xQ'\).  Multiplying by \(x^k\) and summing gives the key identity
\begin{equation}\label{eq:Qrec}
 D^s v=x(v-nQ).
\end{equation}
The coefficient of \(x^{n+1}\) on the right-hand side
is \(nc_n-nc_n=0\), since the coefficient of \(x^n\)
in \(DQ\) is \(nc_n\).

Suppose now that \(Q(r)=Q'(r)=0\) for some \(r>0\).
Then \(v(0)=v(r)=0\), so Lemma~\ref{lem:dil} yields
\[
 \int_0^r vD^s v\,\frac{{\rm d}x}{x}\le0.
\]
On the other hand, using~\eqref{eq:Qrec} and integrating
by parts, we obtain
\begin{equation}\label{eq:prop4-proof1}
 \int_0^r vD^s v\,\frac{{\rm d}x}{x}
 =\int_0^r\bigl(v^2-nxQQ'\bigr)\,{\rm d}x
 =\int_0^r\left(v^2+\frac n2Q^2\right)\,{\rm d}x.
\end{equation}
Here the second equality is a consequence of $\int_0^r xQQ'\,{\rm d}x=-\frac{1}{2}\int_0^r Q^2{\rm d}x$, where the last equality follows by integration by parts, using \(Q(r)=0\). Using $Q(0)=1$ we see that the right-hand side of~\eqref{eq:prop4-proof1} is strictly positive. The obtained contradiction shows that no positive zero of $Q$ can be multiple.

\medskip
Recall that \(Q_{n,1}\) coincides with the Laguerre polynomial
\(L_n\) and hence has \(n\) simple positive zeros.
We shall now show that \(Q_{n,s}\) has \(n\) simple positive
zeros for every \(s\in[1,2]\).

\smallskip
Let $A\subset[1,2]$ be the set of parameters $s\in [1,2]$ for which $Q_{n,s}$
has $n$ simple positive zeros.
The set $A$ is  relatively open in $[1,2]$: around the
zeros of $Q_{n,s_0}$ at any $s_0\in A$, choose $n$ disjoint intervals contained in
$(0,\infty)$ on whose endpoints $Q_{n,s_0}$ has opposite signs.
By continuity of the coefficients, these signs persist for $s$
sufficiently close to $s_0$. The intermediate value theorem then gives
$n$ distinct positive zeros of $Q_{n,s}$, which are necessarily simple
since its degree is $n$.

\smallskip
To show that \(A\) is closed, let
\((s_j)\) be a sequence in \(A\) converging to \(s_*\in[1,2]\). Denote the zeros of $Q_{n,s_j}$ by
$0<r_{j,1}<\cdots<r_{j,n}$. By Vieta's formulas and the explicit expression~\eqref{eq:Q},
we have
\[
 \sum_{k=1}^n r_{j,k}
 =-\frac{[x^{n-1}]Q_{n,s_j}(x)}{[x^n]Q_{n,s_j}(x)}
 =n^{s_j+1}\le n^3,
\]
and
\[
 \sum_{k=1}^n\frac{1}{r_{j,k}}
 =-\frac{[x^1]Q_{n,s_j}(x)}{[x^0]Q_{n,s_j}(x)}
 =n.
\]
Thus $1/n\le r_{j,k}\le n^3$. Passing to a subsequence, we may assume that
$r_{j,k}\to r_k\in[1/n,n^3]$ as $j\to\infty$,
for every $k=1,\ldots,n$. Since \(Q_{n,s_j}(0)=1\) for every \(j\), we have
\[
 Q_{n,s_j}(x)=\prod_{k=1}^n
 \left(1-\frac{x}{r_{j,k}}\right)
 \to
 \prod_{k=1}^n\left(1-\frac{x}{r_k}\right)
 = Q_{n,s_*}(x), \qquad j\to\infty.
\]
The numbers $r_k\geq 1/n$ are positive. If two of the numbers \(r_1,\ldots,r_n\) coincided, then \(Q_{n,s_*}\) would have
a multiple positive zero, which we have already ruled out. Hence $s_*\in A$, proving that $A$ is closed.

\smallskip
Since \(A\) is nonempty, open and closed in the connected
interval \([1,2]\), we conclude that \(A=[1,2]\).
\end{proof}

We shall now extend Proposition~\ref{prop:Q} to all real
\(s\ge1\). To this end, we use finite free multiplicative
convolution~\cite[Definition~1.4]{MSS}, which, up to a sign convention, coincides with
the classical Schur--Szeg\H{o} composition.
Given two polynomials of degree \(n\),
\[
 A(x)=\sum_{k=0}^n a_kx^k,\qquad
 B(x)=\sum_{k=0}^n b_kx^k,
\]
their \emph{finite free multiplicative convolution} is defined by
\begin{equation}\label{eq:finite-multiplicative-convolution}
 (A\boxtimes_n B)(x)
 =\sum_{k=0}^n(-1)^{n-k}
 \frac{a_kb_k}{\binom nk}x^k.
\end{equation}

By the classical Schur--Szeg\H{o} theorem, if \(A\) and
\(B\) have only positive zeros, then so does
\(A\boxtimes_n B\); see~\cite{SzegoSS}
and~\cite[Theorem~1.6]{MSS}.

Moreover, if all zeros of \(A\) are simple, then all zeros
of \(A\boxtimes_n B\) are simple as well.
This is a special case of the results of Kostov and
Shapiro~\cite[Proposition~1.4 and Theorem~1.6]{KS}
on the multiplicities of zeros of Schur--Szeg\H{o}
compositions. For completeness, we give a short direct argument.

Let \(\mathcal P_n\) denote the set of real polynomials
of degree \(n\) with only positive zeros.
The interior of \(\mathcal P_n\) in the space of real
polynomials of degree at most \(n\), equipped with the
usual coefficient topology, consists precisely of the
polynomials with \(n\) simple positive zeros.
Since all zeros of \(B\) are positive, Vieta's formulas
imply that \(b_k\ne0\) for every \(k=0,\ldots,n\).
Hence the linear map \(A\mapsto A\boxtimes_n B\)
is invertible and therefore open.
By the Schur--Szeg\H{o} theorem, this map sends
\(\mathcal P_n\) into itself and consequently sends
the interior of \(\mathcal P_n\) into itself.
This proves the assertion.

\begin{proposition}\label{prop:Q-general}
For every integer \(n\ge1\) and every real \(s\ge1\),
the polynomial \(Q_{n,s}\) has \(n\) simple positive zeros.
\end{proposition}

\begin{proof}
By~\eqref{eq:finite-multiplicative-convolution},
\[
 Q_{n,s}\boxtimes_n Q_{n,t}=(-1)^nQ_{n,s+t}
\]
for all real \(s,t\).
Given \(s\ge1\), write \(s=s_0+m\), where
\(s_0\in[1,2]\) and \(m\) is a nonnegative integer.
By Proposition~\ref{prop:Q}, both \(Q_{n,s_0}\) and
\(Q_{n,1}\) have \(n\) simple positive zeros.
Applying finite free multiplicative convolution with
\(Q_{n,1}\) successively \(m\) times and using the
simplicity-preserving property established above,
we conclude that \(Q_{n,s}\) has \(n\) simple positive zeros.
\end{proof}

\section{Multiplier sequences and real-rootedness of \texorpdfstring{$P_{n,\gamma}$}{P\_n,gamma}}\label{sec:transfer}
In this section, we use Proposition~\ref{prop:Q-general}
to establish the multiplier sequence property of
\(((k!)^{-s})_{k\ge0}\) for \(s\ge1\) and to prove that
\(P_{n,\gamma}\) has \(n\) simple negative zeros for
every \(\gamma\ge2\).

A \emph{multiplier sequence} is a sequence
\(\lambda=(\lambda_k)_{k\ge0}\) of real numbers such that
the diagonal operator
\[
 T_\lambda\left(\sum_k c_kx^k\right)
 =\sum_k\lambda_k c_kx^k
\]
maps every real-rooted real polynomial to a real-rooted
polynomial or to the zero polynomial.
If all \(\lambda_k\) are nonzero, the operator \(T_\lambda\)
preserves the degree of a polynomial.
Note also that the termwise product of two multiplier
sequences is again a multiplier sequence, since the
corresponding operators can be applied in succession.

We shall use the following characterization of positive
multiplier sequences.

\begin{theorem}[P\'olya--Schur~\cite{PS}]\label{th:polya-schur}
For a sequence \(\lambda=(\lambda_k)_{k\ge0}\) of positive
real numbers, the following conditions are equivalent.
\begin{enumerate}
\item \(\lambda\) is a multiplier sequence.
\item For every integer \(d\ge1\), the Jensen polynomial
$
 J_d(x)=\sum_{k=0}^d\binom dk\lambda_kx^k
$
has only negative zeros.
\item The exponential generating function $\sum_{k=0}^\infty\lambda_k\frac{z^k}{k!}$
is an entire function belonging to \(\LP\). Here \(\LP\) is the positive part of the
Laguerre--P\'olya class of type~I, consisting of entire functions that are locally
uniform limits of polynomials with nonnegative coefficients
and only real nonpositive zeros.
\end{enumerate}
\end{theorem}

\begin{proposition}\label{prop:multiplier}
For every real \(s\ge1\), the sequence
$\left(1/(k!)^s\right)_{k\ge0}$
is a multiplier sequence.
\end{proposition}

\begin{proof}
Fix \(s\ge1\). For every integer \(d\ge1\), the Jensen
polynomial associated with the sequence
\(((k!)^{-s})_{k\ge0}\) is given by
\[
 \sum_{k=0}^d\binom dk\frac{x^k}{(k!)^s}
 =Q_{d,s}(-x).
\]
By Proposition~\ref{prop:Q-general}, this polynomial has
only negative zeros. The assertion now follows from the
P\'olya--Schur theorem.
\end{proof}
For every real \(s\), define the diagonal operator \(\T_s: \mathbb R[x] \to \mathbb R[x]\) by
\[
 \T_s\left(\sum_k c_kx^k\right)
 =\sum_k\frac{c_k}{(k!)^s}x^k.
\]
By Proposition~\ref{prop:multiplier}, the operator
\(\T_s\) preserves real-rootedness for every \(s\ge1\).
We are now in a position to prove the real-rootedness
of \(P_{n,\gamma}\).

\begin{proposition}\label{prop:real-rootedness}
For every integer \(n\ge1\) and every real \(\gamma\ge2\),
the polynomial \(P_{n,\gamma}\) has \(n\) simple negative zeros.
\end{proposition}

\begin{proof}
Set \(s=\gamma-1\ge1\) and consider the polynomials
\[
 A(x)=Q_{n,s}(x)
 =\sum_{k=0}^n(-1)^k\binom nk\frac{x^k}{(k!)^s}
\]
and
\[
 B(x)=x^nA(1/x)
 =\sum_{k=0}^n(-1)^{n-k}\binom nk
 \frac{x^k}{((n-k)!)^s}.
\]
By Proposition~\ref{prop:Q-general}, \(A\) has \(n\)
simple positive zeros. Since the zeros of \(B\) are
the reciprocals of the zeros of \(A\), the same holds
for \(B\).
By~\eqref{eq:finite-multiplicative-convolution}, we have
\begin{equation}\label{eq:double}
 (A\boxtimes_n B)(x)
 =\sum_{k=0}^n(-1)^k\binom nk
 \frac{x^k}{(k!)^s((n-k)!)^s}
 =(n!)^{-s}P_{n,s+1}(-x).
\end{equation}
The Schur--Szeg\H{o} theorem and the simplicity-preserving
property established above imply that \(P_{n,\gamma}(-x)\)
has \(n\) simple positive zeros. Equivalently,
\(P_{n,\gamma}(x)\) has \(n\) simple negative zeros.
\end{proof}

\section{Strict interlacing}\label{sec:inter}

\begin{proposition}\label{prop:interlacing}
For every integer \(n\ge1\) and every real \(\gamma\ge2\),
the zeros of \(P_{n,\gamma}\) and \(P_{n+1,\gamma}\)
strictly interlace.
\end{proposition}

\begin{proof}
We first recall that, by the same openness argument as
in Section~\ref{sec:lag}, an invertible, degree-preserving
linear operator that preserves real-rootedness also maps
polynomials with simple real zeros to polynomials with
simple real zeros.
In particular, by Proposition~\ref{prop:multiplier},
this applies to the operator \(\T_s\) for every \(s\ge1\),
since all its diagonal entries are nonzero.

Fix \(\gamma\ge2\) and set \(s=\gamma-1\).
For every integer \(d\ge1\), define the reversed Jensen
polynomial
\[
 R_d(x)=x^dQ_{d,s}(-1/x)
 =\sum_{k=0}^d\binom dk\frac{x^k}{((d-k)!)^s}.
\]
By Proposition~\ref{prop:Q-general}, the polynomial
\(R_d\) has \(d\) simple negative zeros.
A direct comparison of coefficients yields
\begin{equation}\label{eq:pencil0}
 R_{d+1}'(x)=(d+1)R_d(x),\qquad
 \T_sR_d(x)=(d!)^{-s}P_{d,\gamma}(x).
\end{equation}

Recall that if a polynomial \(H\) of degree \(d\) has \(d\)
simple real zeros, then so does \(H+cH'\) for every
\(c\in\RR\); see~\cite[p.~44, Problem~62]{PolyaSzegoII}.
Together with the first identity in~\eqref{eq:pencil0},
this implies that \(R_{n+1}+cR_n\) has \(n+1\) simple
real zeros for every \(c\in\RR\).

Using the second identity in~\eqref{eq:pencil0}, we obtain
\[
 P_{n+1,\gamma}+cP_{n,\gamma}
 =((n+1)!)^s\T_s\left(
 R_{n+1}+\frac{c}{(n+1)^s}R_n\right).
\]
Since \(\T_s\) preserves simple real-rootedness, we conclude that
\begin{equation}\label{eq:pencil}
 P_{n+1,\gamma}+cP_{n,\gamma}
 \text{ has \(n+1\) simple real zeros for every \(c\in\RR\)}.
\end{equation}

By the strict version of the Hermite--Kakeya--Obreschkoff
theorem (see~\cite[Theorem 8]{Branden2004}),
\eqref{eq:pencil} implies that the zeros of \(P_{n,\gamma}\)
and \(P_{n+1,\gamma}\) strictly interlace.
For completeness, we give a short direct argument.
Write \(A=P_{n+1,\gamma}\) and \(B=P_{n,\gamma}\).
First, \(A\) and \(B\) have no common zero.
Indeed, if \(A(r)=B(r)=0\), then \(B'(r)\ne0\)
by Proposition~\ref{prop:real-rootedness}, and choosing
\(c=-A'(r)/B'(r)\) would make \(r\) a multiple zero
of \(A+cB\), contradicting~\eqref{eq:pencil}.

Let \(b_1<\cdots<b_n\) be the zeros of \(B\).
Suppose that two zeros \(u<v\) of \(A\) lie in the same
one of the \(n+1\) intervals
\[
 (-\infty,b_1),\quad (b_1,b_2),\quad\ldots,\quad
 (b_n,\infty).
\]
By Rolle's theorem applied to \(A/B\), there exists
\(\xi\in(u,v)\) such that \((A/B)'(\xi)=0\).
Setting \(c=-A(\xi)/B(\xi)\) and writing
\[
 A+cB=B\left(\frac{A}{B}+c\right),
\]
we see that \(\xi\) is a multiple zero of \(A+cB\),
again contradicting~\eqref{eq:pencil}.

Thus each of the above intervals contains at most one
zero of \(A\). Since \(A\) has \(n+1\) distinct real
zeros and \(A\) and \(B\) have no common zero, each
interval contains exactly one zero of \(A\).
This proves strict interlacing.
\end{proof}

\begin{remark}
The preceding argument also establishes the strict interlacing
of the zeros of \(Q_{n,s}\) and \(Q_{n+1,s}\) for every
\(s\ge1\). Indeed, by the first identity in~\eqref{eq:pencil0}
and Rolle's theorem, the zeros of \(R_n\) and \(R_{n+1}\)
strictly interlace. Since the zeros of \(Q_{d,s}\) are obtained
from those of \(R_d\) by the strictly increasing transformation
\(x\mapsto-1/x\) on \((-\infty,0)\), the zeros of \(Q_{n,s}\)
and \(Q_{n+1,s}\) strictly interlace as well.
\end{remark}

\section{The Le Roy function and its zeros}\label{sec:LeRoy}

We first relate the Le Roy function to the polynomials
\(P_{n,\gamma}\) and deduce the location of its zeros for
\(\gamma\ge2\). We then establish their infinitude and simplicity,
describe the zeros for \(0<\gamma<2\), and obtain sharp criteria
for Laguerre--P\'olya membership and the multiplier sequence property.

\subsection{Polynomial approximation and real-rootedness}

\begin{proposition}[Polynomial approximation]\label{prop:limit}
For every real \(\gamma>0\),
\begin{equation}\label{eq:limit}
 P_{n,\gamma}(z/n^\gamma)\longrightarrow\F_\gamma(z)
\end{equation}
locally uniformly in \(\CC\) as \(n\to\infty\).
\end{proposition}

\begin{proof}
Let \(a_{n,k}\) be the coefficient of \(z^k\) in
\(P_{n,\gamma}(z/n^\gamma)\), with \(a_{n,k}=0\) for \(k>n\).
For \(0\le k\le n\),
\[
 a_{n,k}=\frac{1}{(k!)^\gamma}
 \prod_{j=0}^{k-1}\left(1-\frac{j}{n}\right)^\gamma,
\]
where the empty product for \(k=0\) equals \(1\).
Thus \(a_{n,k}\to1/(k!)^\gamma\) for each fixed \(k\), and
\(0\le a_{n,k}\le1/(k!)^\gamma\) for all \(n,k\).
For every \(R>0\), it follows that
\[
 \sup_{|z|\le R}
 \left|P_{n,\gamma}(z/n^\gamma)-\F_\gamma(z)\right|
 \le\sum_{k=0}^\infty
 \left|a_{n,k}-\frac{1}{(k!)^\gamma}\right|R^k
 \longrightarrow0.
\]
The last convergence follows by dominated convergence, since
\(\sum_{k=0}^\infty R^k/(k!)^\gamma<\infty\).
\end{proof}

\begin{proposition}[Equivalent real-rootedness properties]
\label{prop:equivalence}
For every real \(\gamma>0\), the following conditions are equivalent:
\begin{enumerate}
\item \(P_{n,\gamma}\) is real-rooted for every integer \(n\ge1\).
\item \(\F_\gamma\in\LP\).
\item The sequence \(\bigl((k!)^{1-\gamma}\bigr)_{k\ge0}\)
is a multiplier sequence.
\end{enumerate}
In particular, for every real \(\gamma\ge2\), the function
\(\F_\gamma\) belongs to \(\LP\) and all its zeros are negative.
\end{proposition}

\begin{proof}
Since \(P_{n,\gamma}\) has positive coefficients, condition~1
implies that all its zeros are negative. Proposition~\ref{prop:limit}
and the definition of \(\LP\) then give condition~2.

The equivalence of conditions~2 and~3 follows from the
P\'olya--Schur theorem and the identity
\[
 \F_\gamma(z)=\sum_{k=0}^\infty
 (k!)^{1-\gamma}\frac{z^k}{k!}.
\]

Suppose now that condition~3 holds and put \(s=\gamma-1\).
Then \(\T_s\) preserves real-rootedness. For an integer \(n\ge1\),
set
\[
 A(x)=\T_s(1+x)^n,
 \qquad B(x)=x^nA(1/x).
\]
The polynomial \(A\) is real-rooted and has positive coefficients.
Its reversal \(B\) is therefore also real-rooted. Comparing
coefficients, as in~\eqref{eq:double}, gives
\[
 \T_s B(x)=\sum_{k=0}^n
 \frac{\binom nk}{(k!)^s((n-k)!)^s}x^k
 =(n!)^{-s}P_{n,\gamma}(x).
\]
Thus \(P_{n,\gamma}\) is real-rooted, proving condition~1.

For \(\gamma\ge2\), condition~3 holds by
Proposition~\ref{prop:multiplier}, so \(\F_\gamma\in\LP\).
By Hurwitz's theorem, a nonzero function in \(\LP\) has only
nonpositive zeros. Since \(\F_\gamma(0)=1\), all its zeros
are strictly negative.
\end{proof}

\subsection{Growth and infinitude of zeros}

\begin{proposition}[Growth and factorization]\label{prop:growth}
For every real \(\gamma>0\), the function \(\F_\gamma\) has
order at most \(1/\gamma\). For \(\gamma>1\), it satisfies
\begin{equation}\label{eq:growth}
 \max_{|z|\le R}|\F_\gamma(z)|
 \le e^{\gamma R^{1/\gamma}},\qquad R\ge0,
\end{equation}
has infinitely many zeros, and admits the genus-zero factorization
\[
 \F_\gamma(z)=\prod_{j=1}^\infty
 \left(1-\frac{z}{\xi_j}\right),
 \qquad
 \sum_{j=1}^\infty\frac{1}{|\xi_j|}<\infty,
\]
where the zeros \(\xi_j\) are counted with multiplicities.
\end{proposition}

\begin{proof}
First let \(\gamma>0\) and set \(t=R^{1/\gamma}\).
Since
\[
 \frac{t^k}{k!}=2^{-k}\frac{(2t)^k}{k!}
 \le2^{-k}e^{2t},
\]
we obtain
\[
 \max_{|z|\le R}|\F_\gamma(z)|
 \le\sum_{k=0}^\infty\left(\frac{t^k}{k!}\right)^\gamma
 \le\frac{e^{2\gamma R^{1/\gamma}}}{1-2^{-\gamma}}.
\]
This proves that the order is at most \(1/\gamma\).

For \(\gamma>1\), the inequality
\(\sum_k u_k^\gamma\le(\sum_k u_k)^\gamma\), valid for
nonnegative \(u_k\), gives the sharper estimate
\[
 \max_{|z|\le R}|\F_\gamma(z)|
 \le\sum_{k=0}^\infty\frac{R^k}{(k!)^\gamma}\\
 =\sum_{k=0}^\infty
 \left(\frac{R^{k/\gamma}}{k!}\right)^\gamma
 \le\left(\sum_{k=0}^\infty
 \frac{R^{k/\gamma}}{k!}\right)^\gamma
 =e^{\gamma R^{1/\gamma}}.
\]
In this range the order is less than one. Hadamard's
factorization theorem, together with \(\F_\gamma(0)=1\),
therefore gives the asserted product representation;
see~\cite{Boas}. The product has infinitely many factors,
since \(\F_\gamma\) is not a polynomial.
\end{proof}

For \(\gamma>1\), this factorization also shows that
\(\F_\gamma\in\LP\) if all its zeros are negative: the finite
partial products have nonnegative coefficients and only negative
zeros, and converge locally uniformly to \(\F_\gamma\).
The converse follows from Hurwitz's theorem and \(\F_\gamma(0)=1\).
Thus, for \(\gamma>1\), the conditions in
Proposition~\ref{prop:equivalence} are also equivalent to
the assertion that all zeros of \(\F_\gamma\) are negative.

\subsection{Simplicity of zeros}

Locally uniform convergence does not by itself preserve simplicity
of zeros. We therefore establish simplicity separately.
For \(s\ge1\), extend \(\T_s\) to entire functions by
\[
 \T_s f(z)=\sum_{k=0}^\infty\frac{c_k}{(k!)^s}z^k,
 \qquad f(z)=\sum_{k=0}^\infty c_kz^k.
\]
The resulting series is entire. We shall use the following
extension of the polynomial simplicity-preservation argument
from Section~\ref{sec:inter}.

\begin{proposition}[Preservation of simplicity]
\label{prop:entire-simplicity}
Let \(s\ge1\), and let \(f\) be a real entire function of
order less than one with \(f(0)\ne0\). If all zeros of \(f\)
are simple and negative, then all zeros of \(\T_s f\)
are simple and negative.
\end{proposition}

\begin{proof}
Dividing \(f\) by \(f(0)\), we may assume that \(f(0)=1\).
If \(f\) has finitely many zeros, Hadamard's factorization
theorem makes \(f\) a polynomial, and the assertion follows
from the polynomial case in Section~\ref{sec:inter}.
This includes the case of a nonzero constant.

Otherwise, Hadamard's factorization theorem gives
\[
 f(z)=\prod_{j=1}^\infty\left(1+\frac{z}{\rho_j}\right),
 \qquad \sum_{j=1}^\infty\frac1{\rho_j}<\infty,
\]
where the \(\rho_j\) are distinct positive numbers. Set
\[
 p_N(z)=\prod_{j=1}^N\left(1+\frac{z}{\rho_j}\right),
 \qquad
 F_N(z)=\prod_{j>N}\left(1+\frac{z}{\rho_j}\right).
\]
Then \(f=p_NF_N\), \(p_N\to f\), and \(F_N\to1\),
locally uniformly in \(\CC\).

The operator \(\T_s\) is continuous for locally uniform
convergence. Indeed, writing
\(M_R(h)=\max_{|z|\le R}|h(z)|\), Cauchy's coefficient
estimate gives
\[
 M_R(\T_s h)
 \le M_R(h)\sum_{k=0}^\infty\frac1{(k!)^s}
 \le eM_R(h),\qquad R>0.
\]
Consequently,
\[
 \T_s p_N\longrightarrow g:=\T_s f,
 \qquad H_N:=\T_sF_N\longrightarrow1
\]
locally uniformly. By Proposition~\ref{prop:multiplier} and
Hurwitz's theorem, all zeros of \(g\) are real. They are negative
because \(g\) has nonnegative Taylor coefficients and \(g(0)=1\).

It remains to exclude multiple zeros of \(g\).
For each fixed \(N\), the polynomial \(p_N+t\) has only
simple negative zeros for all real \(t\) sufficiently close
to zero. Restricting also to \(|t|<1/2\), we claim that
\[
 \T_s\bigl((p_N+t)F_N\bigr)=g+tH_N
\]
has only real zeros. To see this, approximate \(F_N\) by
its finite partial products, multiply by \(p_N+t\), and
apply \(\T_s\). The resulting polynomials are real-rooted.
Continuity of \(\T_s\) and Hurwitz's theorem give the claim,
since \((g+tH_N)(0)=1+t\ne0\).

Suppose that \(g\) has a zero \(r\) of multiplicity \(q\ge2\).
Choose \(N\) sufficiently large that \(b:=H_N(r)\ne0\), and write
\[
 g(z)=a(z-r)^q+O\bigl((z-r)^{q+1}\bigr),
 \qquad a\in\RR\setminus\{0\}.
\]
For \(\varepsilon>0\), put
\(t_\varepsilon=(a/b)\varepsilon^q\). Then
\[
 \frac{(g+t_\varepsilon H_N)(r+\varepsilon w)}
      {a\varepsilon^q}
 \longrightarrow w^q+1
\]
locally uniformly in \(w\) as \(\varepsilon\downarrow0\).
For all sufficiently small \(\varepsilon\), the functions
on the left have only real zeros. Since \(w^q+1\) has a
nonreal zero for every \(q\ge2\), this contradicts Hurwitz's
theorem. Hence all zeros of \(g=\T_s f\) are simple.
\end{proof}

\begin{proposition}[Simplicity of the Le Roy zeros]
\label{prop:LeRoy-simplicity}
For every real \(\gamma\ge2\), all zeros of \(\F_\gamma\)
are simple.
\end{proposition}

\begin{proof}
By Proposition~\ref{prop:equivalence}, all zeros of
\(\F_\gamma\) are negative. We distinguish two cases.

\medskip
\noindent\emph{Case 1: \(2\le\gamma\le3\).}
Set \(y(x)=\F_\gamma(-x)\) and \(v=Dy\).
Comparison of coefficients gives
\[
 D^\gamma y=-xy,
 \qquad D^{\gamma-1}v=-xy.
\]
If \(y\) had a multiple zero at \(r>0\), then
\(y(r)=y'(r)=0\), and hence \(v(0)=v(r)=0\).
Since \(1\le\gamma-1\le2\), Lemma~\ref{lem:dil} would yield
\[
 \int_0^r vD^{\gamma-1}v\,\frac{{\rm d}x}{x}\le0.
\]
On the other hand,
\[
 \int_0^r vD^{\gamma-1}v\,\frac{{\rm d}x}{x}
 =-\int_0^r xy'(x)y(x)\,{\rm d}x\\
 =\frac12\int_0^r y(x)^2\,{\rm d}x>0.
\]
The second equality follows by integration by parts using
\(y(r)=0\), and the strict inequality follows from \(y(0)=1\).
This contradiction excludes multiple zeros.

\medskip
\noindent\emph{Case 2: \(\gamma>3\).}
The identity
\[
 \F_\gamma=\T_{\gamma-2}\F_2
\]
follows directly from the power series. By
Proposition~\ref{prop:growth}, the function \(\F_2\) has
order at most \(1/2\), and by Case~1 all its zeros are
simple and negative. Since \(\gamma-2>1\),
Proposition~\ref{prop:entire-simplicity} completes the proof.
\end{proof}

Together, Propositions~\ref{prop:equivalence}, \ref{prop:growth},
and~\ref{prop:LeRoy-simplicity} prove that, for every
\(\gamma\ge2\), the function \(\F_\gamma\) belongs to \(\LP\)
and has infinitely many zeros, all simple and negative.

\subsection{Zeros below the threshold}

At \(\gamma=1\), we have \(\F_1(z)=e^z\), which has no zeros.
The remaining positive parameters below \(2\) behave as follows.

\begin{proposition}\label{prop:LeRoy-subcritical}
\begin{enumerate}
\item For every real \(0<\gamma<1\), the function \(\F_\gamma\)
is positive on \(\RR\) and has infinitely many zeros,
all nonreal.
\item For every real \(1<\gamma<2\), the function \(\F_\gamma\)
has a finite, nonzero number of real zeros, all negative,
and infinitely many nonreal zeros. The number of real zeros,
counted with multiplicities, is odd.
\end{enumerate}
\end{proposition}

\begin{proof}
\emph{Case 1: \(0<\gamma<1\).}
We first prove positivity on the negative real axis using
a variant of the dilation identity from Section~\ref{sec:dilation}.
Set \(y(x)=\F_\gamma(-x)\) and
\[
 K_\gamma=\int_0^\infty(1-e^{-h})
 \frac{{\rm d}h}{h^{1+\gamma}}.
\]
This constant is positive and finite. The substitution \(u=kh\)
gives
\[
 \int_0^\infty(1-e^{-kh})\frac{{\rm d}h}{h^{1+\gamma}}
 =k^\gamma K_\gamma,\qquad k\ge1.
\]
Applying this identity to the power series of \(y\), we obtain
\[
 K_\gamma D^\gamma y(x)
 =\int_0^\infty\bigl(y(x)-y(e^{-h}x)\bigr)
 \frac{{\rm d}h}{h^{1+\gamma}}.
\]
Termwise integration is justified by absolute convergence:
for every \(R>0\),
\[
 \sum_{k=1}^\infty\frac{R^k}{(k!)^\gamma}
 \int_0^\infty(1-e^{-kh})\frac{{\rm d}h}{h^{1+\gamma}}
 =K_\gamma\sum_{k=1}^\infty
 \frac{k^\gamma R^k}{(k!)^\gamma}<\infty.
\]

As before, comparison of coefficients gives \(D^\gamma y=-xy\).
If \(y\) had a positive zero, let \(r\) be its first such zero.
Then \(y(x)>0\) for \(0\le x<r\), since \(y(0)=1\).
Evaluating the integral identity at \(r\) would give
\[
 0=K_\gamma D^\gamma y(r)
 =-\int_0^\infty y(e^{-h}r)
 \frac{{\rm d}h}{h^{1+\gamma}}<0,
\]
a contradiction. Thus \(\F_\gamma\) is positive on the negative
real axis. Its positivity on \([0,\infty)\) follows from its
power series, so it has no real zeros. This positivity is also
a consequence of the probabilistic representation in
\cite[Proposition~9]{BS}.

To prove infinitude of zeros, we use the asymptotic formula
of Boudabsa and Simon~\cite[Proposition~8]{BS}:
\begin{equation}\label{eq:LeRoy-asymptotic-small}
 \F_\gamma(-x)\sim
 \frac{1}{\Gamma(1-\gamma)x(\log x)^\gamma},
 \qquad x\to+\infty.
\end{equation}
Suppose that \(\F_\gamma\) had only finitely many zeros.
It has finite order by Proposition~\ref{prop:growth}, so
Hadamard's factorization theorem would give
$
 \F_\gamma(z)=p(z)e^{q(z)}
$
for polynomials \(p\ne0\) and \(q\). Writing \(m=\deg p\),
we would have
\[
 \log|\F_\gamma(-x)|
 =\operatorname{Re}q(-x)+m\log x+O(1).
\]
However, \eqref{eq:LeRoy-asymptotic-small} gives
\[
 \log|\F_\gamma(-x)|
 =-\log x-\gamma\log\log x+O(1).
\]
The real polynomial \(\operatorname{Re}q(-x)\) must therefore
be constant, since no nonconstant polynomial is \(O(\log x)\).
The two expressions are then incompatible because \(m\ge0\).
Hence \(\F_\gamma\) has infinitely many zeros, all nonreal.

\medskip
\noindent\emph{Case 2: \(1<\gamma<2\).}
Boudabsa and Simon~\cite[Section~5, equation~(15)]{BS}
record the negative-axis asymptotic
\begin{equation}\label{eq:LeRoy-asymptotic-middle}
 \F_\gamma(-x)\sim-\frac{c_\gamma}{x(\log x)^\gamma},
 \qquad x\to+\infty,
\end{equation}
where \(c_\gamma>0\).
Thus \(\F_\gamma\) is strictly negative sufficiently far
to the left, while it is positive on \([0,\infty)\).
All its real zeros lie in a bounded interval and are therefore
finite in number, because \(\F_\gamma\) is a nonzero entire
function. There is at least one such zero, and their total
multiplicity is odd: a real zero changes the sign precisely
when its multiplicity is odd, and the signs at the two ends
of the real axis are opposite.

By Proposition~\ref{prop:growth}, the function \(\F_\gamma\)
has infinitely many zeros altogether. Consequently, infinitely
many of them are nonreal.
\end{proof}

\subsection{Sharp criteria}

The preceding results give a complete criterion for membership
in \(\LP\). They also characterize membership in the full
Laguerre--P\'olya class \(\mathcal{LP}\), consisting of locally
uniform limits of real-rooted real polynomials.

\begin{proposition}[Laguerre--P\'olya criterion]
\label{prop:LeRoy-classification}
For every real \(\gamma>0\), the following conditions are equivalent:
\begin{enumerate}
\item \(\F_\gamma\in\LP\).
\item \(\F_\gamma\in\mathcal{LP}\).
\item \(\F_\gamma\) has no nonreal zeros.
\item \(\gamma\in\{1\}\cup[2,\infty)\).
\end{enumerate}
\end{proposition}

\begin{proof}
Condition~1 implies condition~2, which implies condition~3
by Hurwitz's theorem. Proposition~\ref{prop:LeRoy-subcritical}
shows that condition~3 implies condition~4.
Finally, if \(\gamma\ge2\), condition~1 follows from
Proposition~\ref{prop:equivalence}. If \(\gamma=1\), then
\(\F_1(z)=e^z\in\LP\), as follows from
\((1+z/n)^n\to e^z\) locally uniformly. Thus condition~4
implies condition~1.
\end{proof}

\begin{proposition}[Multiplier sequences]
\label{prop:multiplier-classification}
For a real number \(s\), the sequence
$
 \left(1/(k!)^s\right)_{k\ge0}
$
is a multiplier sequence if and only if
\(s\in\{0\}\cup[1,\infty)\).
\end{proposition}

\begin{proof}
For \(s\ge1\), the assertion follows from
Proposition~\ref{prop:multiplier}. For \(s=0\), the
associated diagonal operator is the identity.

If \(0<s<1\), the exponential generating function of
the sequence is
\[
 \sum_{k=0}^\infty\frac{1}{(k!)^s}\frac{z^k}{k!}
 =\F_{s+1}(z).
\]
Since \(1<s+1<2\), Proposition~\ref{prop:LeRoy-classification}
shows that this function does not belong to \(\LP\).
The P\'olya--Schur theorem therefore excludes the multiplier
sequence property.

Finally, if \(s<0\), applying the corresponding diagonal
operator to \((1+x)^2\) gives
$
 1+2x+2^{-s}x^2,
$
whose discriminant \(4(1-2^{-s})\) is negative. Thus the
sequence is not a multiplier sequence in this case either.
\end{proof}

The nonreal zeros in Proposition~\ref{prop:LeRoy-subcritical}
will be used in Section~\ref{sec:sharp}, together with
polynomial approximation, to study the nonreal zeros of
\(P_{n,\gamma}\) for \(1<\gamma<2\) and of \(Q_{n,s}\)
for \(0<s<1\).

\begin{proof}[Proof of Corollary~\ref{cor:PF}]
Suppose first that \(\gamma>0\). The ordinary generating
function of the sequence \(((k!)^{-\gamma})_{k\ge0}\)
is the entire function
\[
 \F_\gamma(z)=\sum_{k=0}^\infty\frac{z^k}{(k!)^\gamma}.
\]
By the Aissen--Schoenberg--Whitney--Edrei
theorem~\cite{ASW,Edrei}, a sequence of nonnegative real
numbers with an entire generating function is a P\'olya
frequency sequence if and only if its generating function
belongs to \(\LP\).
By Proposition~\ref{prop:LeRoy-classification}, this is
the case precisely when \(\gamma\in\{1\}\cup[2,\infty)\).

For \(\gamma=0\), the sequence \((1)_{k\ge0}\) is a
P\'olya frequency sequence by the same theorem, since
its generating function is \(1/(1-z)\).
Finally, if \(\gamma<0\), the minor of the corresponding
infinite Toeplitz matrix formed by rows \(0,1\) and
columns \(1,2\) is
\[
 \begin{vmatrix}
  1 & 2^{-\gamma}\\
  1 & 1
 \end{vmatrix}
 =1-2^{-\gamma}<0.
\]
Hence the sequence \(((k!)^{-\gamma})_{k\ge0}\)
is not a P\'olya frequency sequence. This completes the proof.
\end{proof}

\section{Zeros below the real-rootedness thresholds}\label{sec:sharp}

We now examine the zeros of \(P_{n,\gamma}\) for \(0<\gamma<2\)
and of \(Q_{n,s}\) for \(s<1\). A representation by moments of
positive random variables gives the exact number of real zeros when
\(0<\gamma<1\) and \(s<0\). In the remaining intervals
\(1<\gamma<2\) and \(0<s<1\), we use polynomial approximation
of the Le Roy function to show that the number of nonreal zeros tends
to infinity with the degree.

\subsection{Factorial moment sequences}

The following consequence of Berg's theorem provides the moment
representations used below.

\begin{lemma}[Factorial moment sequences]\label{lem:factorial-moments}
For every real \(t>0\), there exists a nonconstant random variable
\(X_t>0\) almost surely such that
\begin{equation}\label{eq:factorial-moments}
 \mathbb E[X_t^k]=(k!)^t
 \qquad\text{for every integer }k\ge0.
\end{equation}
\end{lemma}

\begin{proof}
Berg's theorem \cite[Theorem~1.8]{Berg2005}, applied to the
Bernstein function \(f(u)=u\), gives a probability measure on
\((0,\infty)\) with these moments. If \(X_t\) has this distribution,
then
\[
 \operatorname{Var}(X_t)
 =\mathbb E[X_t^2]-(\mathbb E[X_t])^2
 =2^t-1>0,
\]
so \(X_t\) is nonconstant.
\end{proof}

\subsection{Exact number of real zeros}

\begin{proposition}\label{prop:P-small}
Fix a real \(\gamma\in(0,1)\) and an integer \(n\ge1\).
If \(n\) is even, then \(P_{n,\gamma}(x)>0\) for every
\(x\in\RR\). If \(n\) is odd, then \(P_{n,\gamma}\) has
exactly one real zero, namely \(-1\), and this zero is simple.
\end{proposition}

\begin{proof}
Set \(t=1-\gamma>0\), and let \(X_t,Y_t\) be independent random
variables with the distribution supplied by
Lemma~\ref{lem:factorial-moments}. Expanding the power and using
independence gives
\begin{equation}\label{eq:P-moment}
 \frac{1}{(n!)^t}\mathbb E[(Y_t+xX_t)^n]
 =\sum_{k=0}^n\binom nk
   \frac{(k!)^t((n-k)!)^t}{(n!)^t}\,x^k
 =P_{n,\gamma}(x).
\end{equation}
For every real \(x\), independence also yields
\[
 \operatorname{Var}(Y_t+xX_t)
 =(1+x^2)(2^t-1)>0.
\]
Thus \(Y_t+xX_t\) is not almost surely zero. If \(n\) is even,
its \(n\)-th power is nonnegative and is positive with positive
probability. Equation~\eqref{eq:P-moment} therefore gives
\(P_{n,\gamma}(x)>0\).

If \(n\) is odd, differentiating the polynomial identity
\eqref{eq:P-moment} gives
\[
 P'_{n,\gamma}(x)
 =\frac{n}{(n!)^t}
   \mathbb E\!\left[X_t(Y_t+xX_t)^{n-1}\right]>0.
\]
The strict inequality follows from \(X_t>0\) almost surely and the
fact that \(Y_t+xX_t\) is not almost surely zero. Moreover,
\(Y_t-X_t\) has a symmetric
distribution, so its odd moments vanish. Hence
\(P_{n,\gamma}(-1)=0\), and the strict positivity of the derivative
shows that \(-1\) is the unique real zero and is simple.
\end{proof}

\begin{proposition}\label{prop:Q-negative}
Fix a real \(s<0\) and an integer \(n\ge1\).
If \(n\) is even, then \(Q_{n,s}(x)>0\) for every
\(x\in\RR\). If \(n\) is odd, then \(Q_{n,s}\) has exactly
one real zero, which is positive and simple.
\end{proposition}

\begin{proof}
Set \(t=-s>0\), and choose \(X_t\) as in
Lemma~\ref{lem:factorial-moments}. Expanding the power gives
\begin{equation}\label{eq:Q-moment}
 Q_{n,s}(x)=\mathbb E[(1-xX_t)^n].
\end{equation}
Since \(X_t\) is nonconstant, \(1-xX_t\) is not almost surely
zero for any real \(x\). If \(n\) is even,
\eqref{eq:Q-moment} therefore gives \(Q_{n,s}(x)>0\).

If \(n\) is odd, differentiation gives
\[
 Q'_{n,s}(x)
 =-n\,\mathbb E\!\left[X_t(1-xX_t)^{n-1}\right]<0
 \qquad(x\in\RR).
\]
Thus \(Q_{n,s}\) is strictly decreasing on \(\RR\).
Since \(Q_{n,s}(0)=1\) and its leading coefficient is
\(-(n!)^t\), we have \(Q_{n,s}(x)\to-\infty\) as
\(x\to+\infty\). Consequently, \(Q_{n,s}\) has exactly one
real zero, which is positive and simple.
\end{proof}

In particular, in their respective parameter ranges, both families
have exactly
\[
 2\left\lfloor\frac n2\right\rfloor
 =\begin{cases}
 n,&n\text{ even},\\
 n-1,&n\text{ odd},
 \end{cases}
\]
nonreal zeros, counted with multiplicity.

\subsection{Nonreal zeros in the intermediate parameter ranges}

To treat \(0<s<1\), we first record the approximation of the
Le Roy function by the polynomials \(Q_{n,s}\).

\begin{proposition}[Polynomial approximation]\label{prop:Q-limit}
For every real \(s>-1\),
\begin{equation}\label{eq:Q-limit}
 Q_{n,s}\!\left(\frac zn\right)
 \longrightarrow \F_{s+1}(-z)
 \qquad\text{as }n\to\infty,
\end{equation}
locally uniformly for \(z\in\CC\).
\end{proposition}

\begin{proof}
We have
\[
 Q_{n,s}\!\left(\frac zn\right)
 =\sum_{k=0}^n\frac{(-z)^k}{(k!)^{s+1}}
   \prod_{j=0}^{k-1}\left(1-\frac jn\right),
\]
where the empty product equals \(1\). For each fixed \(k\), the
product tends to \(1\) as \(n\to\infty\), and for \(0\le k\le n\)
it lies in \([0,1]\). Extend the coefficients by zero for \(k>n\).
On every disk \(|z|\le R\), the absolute value of the \(k\)-th term
is bounded by \(R^k/(k!)^{s+1}\). Since \(s+1>0\), these bounds
form a summable series. Termwise convergence and this summable bound
imply uniform convergence on \(|z|\le R\), proving
\eqref{eq:Q-limit}.
\end{proof}

\begin{proposition}\label{prop:nonreal}
The following assertions hold.
\begin{enumerate}
\item For every fixed real \(\gamma\in(1,2)\), the number of
nonreal zeros of \(P_{n,\gamma}\) tends to infinity as
\(n\to\infty\).
\item For every fixed real \(s\in(0,1)\), the number of nonreal
zeros of \(Q_{n,s}\) tends to infinity as \(n\to\infty\).
\end{enumerate}
In both assertions, zeros are counted with multiplicity.
\end{proposition}

\begin{proof}
Fix \(\alpha\in(1,2)\), and set
\[
 U_n(z)=P_{n,\alpha}\!\left(\frac{z}{n^\alpha}\right),
 \qquad
 V_n(z)=Q_{n,\alpha-1}\!\left(-\frac zn\right).
\]
By Propositions~\ref{prop:limit} and~\ref{prop:Q-limit}, both
sequences converge locally uniformly to \(\F_\alpha\).
Proposition~\ref{prop:LeRoy-subcritical} shows that
\(\F_\alpha\) has infinitely many nonreal zeros.

Let \(M\ge1\) be an integer. Choose \(M\) distinct nonreal zeros
of \(\F_\alpha\), and surround them by pairwise disjoint closed
disks \(D_1,\ldots,D_M\) contained in \(\CC\setminus\RR\),
with no zeros of \(\F_\alpha\) on their boundaries. By locally
uniform convergence, for all sufficiently large \(n\) and every
\(z\in\bigcup_{j=1}^M\partial D_j\),
\[
 |U_n(z)-\F_\alpha(z)|<|\F_\alpha(z)|,
 \qquad
 |V_n(z)-\F_\alpha(z)|<|\F_\alpha(z)|.
\]
Rouch\'e's theorem implies that \(U_n\) and \(V_n\) each have
the same number of zeros as \(\F_\alpha\) in every disk,
counted with multiplicity. Each polynomial therefore has at least
\(M\) nonreal zeros for all sufficiently large \(n\).

Rescaling the variable by a nonzero real number preserves
nonreality of zeros. Since \(M\) is arbitrary, the first assertion
follows by taking \(\alpha=\gamma\), and the second by taking
\(\alpha=s+1\).
\end{proof}

\parhead{Real-rootedness criteria.}
At $\gamma=1$ and $s=0$, we have
\(P_{n,1}(x)=(1+x)^n\) and  \(Q_{n,0}(x)=(1-x)^n\).
Combining these identities and the preceding results with
Propositions~\ref{prop:real-rootedness} and~\ref{prop:Q-general}
gives, for \(\gamma>0\) and \(s\in\RR\),
\[
\begin{aligned}
 P_{n,\gamma}\text{ is real-rooted for every integer }n\ge1
 &\iff \gamma\in\{1\}\cup[2,\infty),\\
 Q_{n,s}\text{ is real-rooted for every integer }n\ge1
 &\iff s\in\{0\}\cup[1,\infty).
\end{aligned}
\]
The results on zeros are summarized below.
\begin{center}
\begin{tabular}{@{}lp{0.70\linewidth}@{}}
\toprule
\multicolumn{2}{c}{The polynomials \(P_{n,\gamma}\)}\\
\midrule
Parameter & Zeros\\
\midrule
\(0<\gamma<1\)
& No real zeros for even \(n\); exactly one simple real zero,
  at \(-1\), for odd \(n\).\\
\addlinespace
\(\gamma=1\)
& All zeros coincide at \(-1\).\\
\addlinespace
\(1<\gamma<2\)
& The number of nonreal zeros tends to infinity as $n\to\infty$.\\
\addlinespace
\(\gamma\ge2\)
& All \(n\) zeros are simple and negative; the zeros of consecutive degrees
  strictly interlace.\\
\bottomrule
\end{tabular}
\end{center}

\begin{center}
\begin{tabular}{@{}lp{0.70\linewidth}@{}}
\toprule
\multicolumn{2}{c}{The polynomials \(Q_{n,s}\)}\\
\midrule
Parameter & Zeros\\
\midrule
\(s<0\)
& No real zeros for even \(n\); exactly one simple positive
  real zero for odd \(n\).\\
\addlinespace
\(s=0\)
& All zeros coincide at \(1\).\\
\addlinespace
\(0<s<1\)
& The number of nonreal zeros tends to infinity as $n\to\infty$.\\
\addlinespace
\(s\ge1\)
& All \(n\) zeros are simple and positive; the zeros
  of consecutive degrees strictly interlace.\\
\bottomrule
\end{tabular}
\end{center}

\section{Limiting empirical distributions of zeros}\label{sec:positions}

Following~\cite{JKM1,JKM}, for a polynomial \(A\)
of degree \(n\ge1\) with zeros \(z_1,\ldots,z_n\), counted with
multiplicity, we write
\[
 \lsem A\rsem_n=\frac1n\sum_{k=1}^n\delta_{z_k}
\]
for its empirical zero distribution. We use
\(\overset{\mathrm w}{\longrightarrow}\) to denote weak convergence
of probability measures. We first determine the limiting distributions
for \(P_{n,\gamma}\) and \(Q_{n,s}\), and then describe how
factorial multiplier operators transform a general limiting zero
distribution.

\parhead{The exponential-profile method.}
We use the following consequence of the exponential-profile theorem
of Jalowy, Kabluchko, and Marynych \cite[Theorem~2.4]{JKM}.
For every integer \(n\ge1\), let
\[
 A_n(x)=\sum_{k=0}^n a_{n,k}x^k
\]
have degree \(n\), nonnegative coefficients, and only nonpositive
zeros. Suppose that
\[
 \frac1n\log a_{n,\lfloor nu\rfloor}\longrightarrow h(u)
 \qquad(0<u<1),
\]
where \(h\) is finite on \((0,1)\) and \(\log0=-\infty\).
Then \(\lsem A_n\rsem_n\)
converges weakly to a probability measure \(\eta\) on
\((-\infty,0)\). Its Cauchy transform
\[
 G_\eta(t)=\int_{-\infty}^0\frac{{\rm d}\eta(x)}{t-x},
 \qquad t>0,
\]
is characterized by the fact that \(t\mapsto tG_\eta(t)\) is the
inverse of \(u\mapsto e^{-h'(u)}\). The smoothness of \(h\) and
the existence of this inverse are part of the theorem.

\subsection{Limiting zero distribution of \texorpdfstring{\(P_{n,\gamma}\)}{P\_n,gamma}}
The following result extends the integer-exponent case of
\cite[Theorem~4.7]{JKM} to every real \(\gamma\ge2\).

\begin{theorem}[Limit distribution for zeros of $P_{n,\gamma}$]\label{th:P-zero-limit}
For every real \(\gamma\ge2\),
\[
 \lsem P_{n,\gamma}\rsem_n
 \overset{\mathrm w}{\longrightarrow}\mu_\gamma
 \qquad\text{as }n\to\infty,
\]
where \(\mu_\gamma\) is the probability measure on \((-\infty,0)\)
with density
\begin{equation}\label{eq:P-limit-density}
 p_\gamma(x)=
 \frac{\sin(\pi/\gamma)\,|x|^{1/\gamma-1}}
 {\pi\bigl(1+2|x|^{1/\gamma}\cos(\pi/\gamma)
                 +|x|^{2/\gamma}\bigr)},
 \qquad x<0.
\end{equation}
\end{theorem}

\begin{proof}
Proposition~\ref{prop:real-rootedness} supplies the required
negative zeros. By Stirling's formula, the exponential profile is
\[
 h_\gamma(u)
 =-\gamma\bigl(u\log u+(1-u)\log(1-u)\bigr),
 \qquad 0<u<1.
\]
Thus
\[
 e^{-h_\gamma'(u)}=\left(\frac{u}{1-u}\right)^\gamma.
\]
The exponential-profile theorem gives weak convergence and identifies
the Cauchy transform of the limit as
\[
 G_\gamma(t)=\frac{t^{1/\gamma-1}}{1+t^{1/\gamma}},
 \qquad t>0.
\]
This identity extends analytically to
\(\CC\setminus(-\infty,0]\), using the branch of
\(z^{1/\gamma}\) that is positive for positive \(z\).
For \(x<0\), the upper boundary value is obtained by substituting
\(z^{1/\gamma}=|x|^{1/\gamma}e^{i\pi/\gamma}\).
These boundary values exist locally uniformly on \((-\infty,0)\),
so the Stieltjes--Perron inversion formula gives
\[
 p_\gamma(x)=-\frac1\pi\operatorname{Im}G_\gamma(x+i0),
\]
which is \eqref{eq:P-limit-density}.
\end{proof}

\subsection{The modified Laguerre polynomials}

To describe the next limit, for every real \(s>0\) let \(\nu_s\)
be the \emph{Fuss--Catalan distribution} with moments
\begin{equation}\label{eq:FC-moments}
 \int_0^\infty x^k\,{\rm d}\nu_s(x)
 =\frac{1}{sk+1}\binom{(s+1)k}{k},
 \qquad k=0,1,2,\ldots,
\end{equation}
where the binomial coefficient is understood in the generalized
sense. This is a unique compactly supported probability measure.
It is absolutely continuous, with support
\begin{equation}\label{eq:FC-support}
 [0,b_s],\qquad b_s=\frac{(s+1)^{s+1}}{s^s};
\end{equation}
see \cite[Theorem~2.4 and the subsequent remark]{MPZ}.
In particular, \(\nu_1\) is the standard Marchenko--Pastur
distribution \(\mathrm{MP}\), whose density is
\[
 \frac{1}{2\pi}\sqrt{\frac{4-x}{x}},\qquad 0<x<4.
\]
More generally,
\begin{equation}\label{eq:FC-free-power}
 \nu_s=\mathrm{MP}^{\boxtimes s},\qquad s>0,
\end{equation}
where \(\boxtimes\) denotes free multiplicative convolution;
see \cite[Introduction]{MPZ}.

For the proof below, we use the transform
\[
 M_s(z)=\int_0^{b_s}\frac{{\rm d}\nu_s(x)}{1-zx},
 \qquad |z|<b_s^{-1}.
\]
By Lagrange inversion, the moment formula \eqref{eq:FC-moments}
is equivalent to
\begin{equation}\label{eq:FC-transform}
 M_s(z)=1+zM_s(z)^{s+1},\qquad M_s(0)=1,
\end{equation}
in a neighborhood of zero, with the power defined by its analytic
branch near \(M_s(0)=1\).

\begin{theorem}[Limit distribution for zeros of $Q_{n,s}(n^s\,\cdot)$]\label{th:Q-zero-limit}
For every real \(s\ge1\),
\begin{equation}\label{eq:Q-zero-limit}
 \lsem Q_{n,s}(n^s\,\cdot)\rsem_n
 \overset{\mathrm w}{\longrightarrow}\nu_s
 \qquad\text{as }n\to\infty.
\end{equation}
Thus each zero of \(Q_{n,s}\) is divided by \(n^s\) before
forming the empirical distribution.
\end{theorem}

\begin{proof}
Consider the reflected and rescaled polynomials
\[
 R_{n,s}(x)=Q_{n,s}(-n^sx)
 =\sum_{k=0}^n\binom nk\frac{n^{sk}}{(k!)^s}x^k.
\]
They have positive coefficients and only negative zeros by
Proposition~\ref{prop:Q-general}. Stirling's formula gives the
exponential profile
\[
 h_s(u)=-(s+1)u\log u-(1-u)\log(1-u)+su,
 \qquad 0<u<1,
\]
and hence
\[
 e^{-h_s'(u)}=\frac{u^{s+1}}{1-u}.
\]
The exponential-profile theorem gives a weak limit for
\(\lsem R_{n,s}\rsem_n\). Let \(\rho_s\) be its reflection
about zero. The transform
\[
 H_{\rho_s}(t)=\int_0^\infty\frac{t}{t+x}\,{\rm d}\rho_s(x),
 \qquad t>0,
\]
takes values in \((0,1)\) and satisfies
\begin{equation}\label{eq:FC-H-transform}
 H_{\rho_s}(t)^{s+1}=t\bigl(1-H_{\rho_s}(t)\bigr).
\end{equation}

On the other hand, \(H_{\nu_s}(t)=M_s(-1/t)\), so
\eqref{eq:FC-transform} shows that \(H_{\nu_s}\) satisfies the
same equation for all sufficiently large \(t\). The function
\(u\mapsto u^{s+1}/(1-u)\) is strictly increasing on \((0,1)\).
Therefore \(H_{\rho_s}(t)=H_{\nu_s}(t)\) for all sufficiently
large \(t\). Uniqueness of the Stieltjes transform gives
\(\rho_s=\nu_s\), proving \eqref{eq:Q-zero-limit}.
\end{proof}

\parhead{The limiting density.}
The density \(q_s\) of \(\nu_s\) has the following parametric
representation \cite[Proposition~2.1]{FL}. Set
\begin{equation}\label{eq:FC-parametrization}
 x_s(\varphi)=
 \frac{\sin^{s+1}((s+1)\varphi)}
 {\sin\varphi\,\sin^s(s\varphi)},
 \qquad 0<\varphi<\frac{\pi}{s+1}.
\end{equation}
Then
\begin{equation}\label{eq:FC-density}
 q_s(x_s(\varphi))=
 \frac{\sin^2\varphi\,\sin^{s-1}(s\varphi)}
 {\pi\sin^s((s+1)\varphi)}.
\end{equation}
As \(\varphi\) increases from \(0\) to \(\pi/(s+1)\),
\(x_s(\varphi)\) decreases from \(b_s\) to \(0\).

\subsection{Multiplier sequences and finite free convolution}

The preceding limit also describes how factorial multipliers act on
general sequences of polynomials with nonnegative zeros.
For a polynomial \(A(x)=\sum_{k=0}^n a_kx^k\) of degree \(n\)
and a real \(s\ge1\), define
\begin{equation}\label{eq:normalized-factorial-operator}
 (\T_{n,s}A)(x)
 =\sum_{k=0}^n a_k\frac{n^{sk}}{(k!)^s}x^k
 =(\T_sA)(n^sx).
\end{equation}

\begin{theorem}\label{th:factorial-zero-limit}
Fix a real \(s\ge1\). For every integer \(n\ge1\), let \(A_n\)
be a real polynomial of degree \(n\) with only nonnegative zeros.
Suppose that
\[
 \lsem A_n\rsem_n
 \overset{\mathrm w}{\longrightarrow}\mu
 \qquad\text{as }n\to\infty,
\]
where \(\mu\) is a probability measure on \([0,\infty)\).
Then \(\T_{n,s}A_n\) has degree \(n\) and only nonnegative
zeros, and
\begin{equation}\label{eq:factorial-zero-limit}
 \lsem\T_{n,s}A_n\rsem_n
 \overset{\mathrm w}{\longrightarrow}
 \mu\boxtimes\nu_s
 =\mu\boxtimes\mathrm{MP}^{\boxtimes s}
 \qquad\text{as }n\to\infty.
\end{equation}
No compact-support or moment assumption on \(\mu\) is required.
\end{theorem}

\begin{proof}
Proposition~\ref{prop:multiplier} gives real-rootedness of
\(\T_{n,s}A_n\), and its degree remains \(n\). All nonzero
coefficients of \(A_n(-x)\) have the same sign. Positive coefficient
multipliers preserve this property, so \(\T_{n,s}A_n\) has no
negative zeros.

Set
\[
 B_{n,s}(x)=(-1)^nQ_{n,s}(n^sx)
 =\sum_{k=0}^n(-1)^{n-k}\binom nk
   \frac{n^{sk}}{(k!)^s}x^k.
\]
Recall the definition of finite free multiplicative convolution
in~\eqref{eq:finite-multiplicative-convolution}. Comparing coefficients in
\eqref{eq:finite-multiplicative-convolution} gives the exact identity
\begin{equation}\label{eq:factorial-convolution-identity}
 \T_{n,s}A_n=A_n\boxtimes_n B_{n,s}.
\end{equation}
Moreover, Theorem~\ref{th:Q-zero-limit} gives
\(\lsem B_{n,s}\rsem_n\overset{\mathrm w}{\longrightarrow}\nu_s\).

Suppose first that \(\mu\ne\delta_0\). Since neither limiting
measure has mass at infinity and \(\nu_s(\{0\})=0\), the
hypotheses of \cite[Theorem~3.7]{JKM1} are satisfied. That theorem
and \eqref{eq:factorial-convolution-identity} imply
\eqref{eq:factorial-zero-limit}.

It remains to consider \(\mu=\delta_0\). Multiplying \(A_n\)
by a nonzero constant, we may assume that it is monic. Write
\[
 f_n(x)=(-1)^nA_n(-x)=\sum_{k=0}^n c_{n,k}x^k,
 \qquad c_{n,k}\ge0,
\]
and set
\[
 g_n(x)=(\T_{n,s}f_n)(x)
 =\sum_{k=0}^n w_{n,k}c_{n,k}x^k,
 \qquad w_{n,k}=\frac{n^{sk}}{(k!)^s}.
\]
The weights are nondecreasing in \(k\), because
\[
 \frac{w_{n,k+1}}{w_{n,k}}
 =\left(\frac{n}{k+1}\right)^s\ge1
 \qquad(0\le k<n).
\]
Consequently,
\[
 g_n'(1)f_n(1)-f_n'(1)g_n(1)
 =\sum_{0\le j<k\le n}
 (k-j)(w_{n,k}-w_{n,j})c_{n,j}c_{n,k}\ge0.
\]
The assumption \(\lsem A_n\rsem_n\to\delta_0\) implies
\[
 \frac{f_n'(1)}{nf_n(1)}
 =\int_0^\infty\frac{\lsem A_n\rsem_n({\rm d}x)}{1+x}
 \longrightarrow1.
\]
Thus
\[
 1\ge
 \int_0^\infty
 \frac{\lsem\T_{n,s}A_n\rsem_n({\rm d}x)}{1+x}
 =\frac{g_n'(1)}{ng_n(1)}
 \ge\frac{f_n'(1)}{nf_n(1)}
 \longrightarrow1.
\]
For every \(\varepsilon>0\), it follows that
\[
 \lsem\T_{n,s}A_n\rsem_n([\varepsilon,\infty))
 \le\frac{1+\varepsilon}{\varepsilon}
 \left(1-\frac{g_n'(1)}{ng_n(1)}\right)
 \longrightarrow0.
\]
This proves convergence to
\(\delta_0=\delta_0\boxtimes\nu_s\) and completes the proof.
\end{proof}

For polynomials with nonpositive zeros, the corresponding statement
follows by applying Theorem~\ref{th:factorial-zero-limit} to
\(A_n(-x)\) and reflecting the limiting measure about zero.
Taking \(A_n(x)=(x-1)^n\), whose empirical zero distribution is
\(\delta_1\), recovers Theorem~\ref{th:Q-zero-limit}.

\section*{Acknowledgements}
ZK was supported by the DFG under Germany's Excellence Strategy
EXC 2044/2--390685587, Mathematics M\"unster:
\emph{Dynamics-Geometry-Structure} and by the DFG Research
Training Group \emph{Rigorous Analysis of Complex Random Systems}
(RTG 3027, Project Number 524444762).

\section*{Declarations}
\subsection*{AI statement}
The proof of the real-rootedness of \(P_{n,\gamma}\) for
\(\gamma\ge2\) was obtained by ChatGPT (GPT-6.2 Astra).
The present paper developed through an extensive mathematical
dialogue between the authors and ChatGPT, involving the
formulation of further results, the development of proofs,
and the refinement of the exposition.
ChatGPT also assisted with the literature search and the
preparation of the manuscript.
The authors take responsibility for the content of the paper.

\subsection*{Conflict of interest statement}
The authors declare that they have no conflict of interest.

\subsection*{Data availability statement}
We do not analyze or generate any datasets.

\bibliographystyle{amsplain}
\bibliography{references}

\providecommand{\bysame}{\leavevmode\hbox to3em{\hrulefill}\thinspace}
\providecommand{\MR}{\relax\ifhmode\unskip\space\fi MR }
\providecommand{\MRhref}[2]{%
  \href{http://www.ams.org/mathscinet-getitem?mr=#1}{#2}
}
\providecommand{\href}[2]{#2}
\begin{thebibliography}{10}

\bibitem{ASW}
Michael Aissen, I.~J. Schoenberg, and A.~M. Whitney, \emph{On the generating
  functions of totally positive sequences, {I}}, J. Analyse Math. \textbf{2}
  (1952), 93--103, \doi{10.1007/BF02786970}.

\bibitem{Berg2005}
Christian Berg, \emph{On powers of {Stieltjes} moment sequences, {I}}, J.
  Theoret. Probab. \textbf{18} (2005), no.~4, 871--889,
  \doi{10.1007/s10959-005-7530-6}.

\bibitem{Boas}
Ralph~Philip Boas, Jr., \emph{Entire functions}, Pure and Applied Mathematics,
  vol.~5, Academic Press, New York, 1954.

\bibitem{BS}
Lotfi Boudabsa and Thomas Simon, \emph{Some properties of the {Kilbas--Saigo}
  function}, Mathematics \textbf{9} (2021), no.~3, Article no.~217, 24~pp.,
  \doi{10.3390/math9030217}.

\bibitem{Branden2004}
Petter Br{\"a}nd{\'e}n, \emph{On operators on polynomials preserving
  real-rootedness and the {Neggers--Stanley} conjecture}, J. Algebraic Combin.
  \textbf{20} (2004), no.~2, 119--130,
  \doi{10.1023/B:JACO.0000047295.93525.df}.

\bibitem{CravenCsordas}
Thomas Craven and George Csordas, \emph{The {Fox--Wright} functions and
  {Laguerre} multiplier sequences}, J. Math. Anal. Appl. \textbf{314} (2006),
  no.~1, 109--125, \doi{10.1016/j.jmaa.2005.03.058}.

\bibitem{Duran}
Antonio~J. Dur{\'a}n, \emph{{Brenke} polynomials with real zeros and the
  {Riemann} hypothesis}, Constr. Approx. (2026), Published online 1 October
  2026. \doi{10.1007/s00365-026-09787-y}. Preprint \arxiv{2405.18940}.

\bibitem{Edrei}
Albert Edrei, \emph{On the generating functions of totally positive sequences,
  {II}}, J. Analyse Math. \textbf{2} (1952), 104--109,
  \doi{10.1007/BF02786971}.

\bibitem{Fisk}
Steve Fisk, \emph{Polynomials, roots, and interlacing},
  \href{https://arxiv.org/abs/math/0612833}{arXiv:math/0612833}, 2006, Revised
  11 March 2008.

\bibitem{FL}
Peter~J. Forrester and Dang-Zheng Liu, \emph{{Raney} distributions and random
  matrix theory}, J. Stat. Phys. \textbf{158} (2015), no.~5, 1051--1082,
  \doi{10.1007/s10955-014-1150-4}.

\bibitem{Gerhold}
Stefan Gerhold, \emph{Asymptotics for a variant of the {Mittag--Leffler}
  function}, Integral Transforms Spec. Funct. \textbf{23} (2012), no.~6,
  397--403, \doi{10.1080/10652469.2011.596151}.

\bibitem{Hille}
Einar Hille, \emph{Note on some hypergeometric series of higher order}, J.
  London Math. Soc. \textbf{4} (1929), no.~1, 50--54,
  \doi{10.1112/jlms/s1-4.1.50}.

\bibitem{Hurwitz}
Adolf Hurwitz, \emph{{\"U}ber die {Wurzeln} einiger transzendenten
  {Gleichungen}}, Mitt. Math. Ges. Hamburg \textbf{2} (1890), 25--31, Reprinted
  in {Mathematische Werke}, vol.~1, 1932, pp.~299--305. Reprint available at
  \url{https://link.springer.com/chapter/10.1007/978-3-0348-4161-0_17}.

\bibitem{Hutchinson}
J.~I. Hutchinson, \emph{On a remarkable class of entire functions}, Trans.
  Amer. Math. Soc. \textbf{25} (1923), no.~3, 325--332,
  \doi{10.1090/S0002-9947-1923-1501248-1}.

\bibitem{JKM1}
Jonas Jalowy, Zakhar Kabluchko, and Alexander Marynych, \emph{Zeros and
  exponential profiles of polynomials {I}: {Limit} distributions, finite free
  convolutions and repeated differentiation}, Constr. Approx. (2026), Published
  online 6 October 2026. \doi{10.1007/s00365-026-09778-z}. Preprint
  \arxiv{2504.11593}.

\bibitem{JKM}
\bysame, \emph{Zeros and exponential profiles of polynomials {II}: {Examples}},
  Constr. Approx. (2026), Published online 11 September 2026.
  \doi{10.1007/s00365-026-09775-2}. Preprint \arxiv{2509.11248}.

\bibitem{KLV}
O.~M. Katkova, T.~Lobova, and A.~M. Vishnyakova, \emph{On power series having
  sections with only real zeros}, Comput. Methods Funct. Theory \textbf{3}
  (2003), no.~2, 425--441, \doi{10.1007/BF03321047}.

\bibitem{KV2024}
Olga Katkova and Anna Vishnyakova, \emph{An analog of multiplier sequences for
  the set of totally positive sequences}, J. Math. Phys. Anal. Geom.
  \textbf{20} (2024), no.~3, 353--371, \doi{10.15407/mag20.03.06}.

\bibitem{Khare}
Apoorva Khare, \emph{Matrix analysis, {P{\'o}lya} frequency functions, and
  preservers of (total) positivity}, Lecture notes, Indian Institute of
  Science, \url{https://math.iisc.ac.in/~khare/teaching/Math341-notes.pdf},
  accessed 29 November 2026.

\bibitem{KS}
Vladimir Kostov and Boris Shapiro, \emph{On the {Schur--Szeg{\"o}} composition
  of polynomials}, C. R. Math. Acad. Sci. Paris \textbf{343} (2006), no.~2,
  81--86, \doi{10.1016/j.crma.2006.06.007}.

\bibitem{LeRoy}
{\'E}douard Le~Roy, \emph{Valeurs asymptotiques de certaines s{\'e}ries
  proc{\'e}dant suivant les puissances enti{\`e}res et positives d'une variable
  r{\'e}elle}, Bull. Sci. Math. (2) \textbf{24} (1900), 245--268.

\bibitem{MSS}
Adam~W. Marcus, Daniel~A. Spielman, and Nikhil Srivastava, \emph{Finite free
  convolutions of polynomials}, Probab. Theory Related Fields \textbf{182}
  (2022), no.~3--4, 807--848, \doi{10.1007/s00440-021-01105-w}.

\bibitem{MPZ}
Wojciech M{\l}otkowski, Karol~A. Penson, and Karol {\.Z}yczkowski,
  \emph{Densities of the {Raney} distributions}, Doc. Math. \textbf{18} (2013),
  1573--1596, \doi{10.4171/DM/437}.

\bibitem{NV2018}
Thu~Hien Nguyen and Anna Vishnyakova, \emph{On the entire functions from the
  {Laguerre--P\'olya} class having the decreasing second quotients of {Taylor}
  coefficients}, J. Math. Anal. Appl. \textbf{465} (2018), no.~1, 348--358,
  \doi{10.1016/j.jmaa.2018.05.018}.

\bibitem{NV2024}
\bysame, \emph{{Hutchinson}'s intervals and entire functions from the
  {Laguerre--P\'olya} class}, Eur. J. Math. \textbf{10} (2024), Article no.~10,
  14~pp., \doi{10.1007/s40879-023-00723-z}.

\bibitem{PS}
George P{\'o}lya and Issai Schur, \emph{{\"Uber zwei Arten von Faktorenfolgen
  in der Theorie der algebraischen Gleichungen}}, J. Reine Angew. Math.
  \textbf{144} (1914), 89--113, \doi{10.1515/crll.1914.144.89}.

\bibitem{PolyaSzegoII}
George P{\'o}lya and G{\'a}bor Szeg{\H{o}}, \emph{Problems and theorems in
  analysis {II}: {Theory} of functions. {Zeros}. {Polynomials}. {Determinants}.
  {Number} theory. {Geometry}}, Classics in Mathematics, Springer-Verlag,
  Berlin, 1998, \doi{10.1007/978-3-642-61905-2}.

\bibitem{Sokal}
Alan~D. Sokal, \emph{When does a hypergeometric function {${}_pF_q$} belong to
  the {Laguerre--P{\'o}lya} class {$LP^+$}?}, J. Math. Anal. Appl. \textbf{515}
  (2022), no.~2, Article no.~126432, 8~pp., \doi{10.1016/j.jmaa.2022.126432}.

\bibitem{StanleyMO}
Richard~P. Stanley, \emph{Is {$(n!^{-d})_{n\geq0}$} a {P{\'o}lya} frequency
  sequence?}, MathOverflow, \url{https://mathoverflow.net/questions/477780},
  2024, Question posted 29 August 2024; with an answer by Mateusz
  Kwa{\'s}nicki, 2 September 2024.

\bibitem{SzegoSS}
G{\'a}bor Szeg{\H{o}}, \emph{{Bemerkungen zu einem Satz von J.~H. Grace
  {\"u}ber die Wurzeln algebraischer Gleichungen}}, Math. Z. \textbf{13}
  (1922), 28--55, \doi{10.1007/BF01485280}.

\bibitem{SzegoOP}
\bysame, \emph{Orthogonal polynomials}, fourth ed., American Mathematical
  Society Colloquium Publications, vol.~23, American Mathematical Society,
  Providence, RI, 1975.

\bibitem{WZ}
Yi~Wang and Bin Zhang, \emph{{Hadamard} powers of polynomials with only real
  zeros}, Linear Algebra Appl. \textbf{439} (2013), no.~10, 3173--3176,
  \doi{10.1016/j.laa.2013.08.012}.

\end{thebibliography}

\end{document}